\documentclass[11pt]{amsart}
\usepackage{amsfonts,amscd,amssymb}
\usepackage[dvips]{epsfig,graphics}

\newtheorem{thm}{Theorem}[section]
\newtheorem{cor}[thm]{Corollary}
\newtheorem{lemma}[thm]{Lemma}

\theoremstyle{definition}
\newtheorem{defn}[thm]{Definition}
\newtheorem{ex}[thm]{Example}

\newcommand{\R}{\mathbb R}

\renewcommand{\o}{\operatorname}

\newcommand{\oto}{\o{O}(2,1)}
\newcommand{\soto}{{\o{SO}(2,1)^0}}

\newcommand{\rto}{\R^{2,1}}

\newcommand{\affine}{{\mathbb A}^{2,1}}
\newcommand{\isoaff}{\o{Aff}(\affine)}
\newcommand{\affn}{{\mathbb A}^n}
\newcommand{\E}{\affn}

\newcommand{\vx}{{\mathsf x}}
\newcommand{\vy}{{\mathsf y}}
\newcommand{\vz}{{\mathsf z}}
\newcommand{\vu}{{\mathsf u}}
\newcommand{\vv}{{\mathsf v}}
\newcommand{\vw}{{\mathsf w}}

\newcommand{\xo}[1]{{\mathsf x}^0_{#1}}
\newcommand{\xp}[1]{{\mathsf x}^+_{#1}}
\newcommand{\xm}[1]{{\mathsf x}^-_{#1}}
\newcommand{\xpm}[1]{{\mathsf x}^{\pm}_{#1}}

\newcommand{\ldot}[1]{\mathbb B(#1)}

\newcommand{\cqfd}{\begin{flushright}$\square$\end{flushright}}

\renewcommand{\S}{{\mathbb S}^{2,1}}

\newcommand{\F}{\mathfrak F}
\newcommand{\PO}[1]{\mathfrak O[#1]}
\newcommand{\UP}{U^+}
\newcommand{\UM}{U^-}
\newcommand{\UPN}{U^+_n}
\newcommand{\C}[1]{{\mathcal C}[#1]}

\newcommand{\HN}{{\mathcal H}_n}

\newcommand{\vg}{\vv_g}
\newcommand{\vh}{\vv_h}
\newcommand{\lag}{\lambda_g}
\newcommand{\lah}{\lambda_h}

\begin{document}

\title{Non-proper Actions of the Fundamental Group of a Punctured Torus}
\author{Virginie Charette}
    \address{Department of Mathematics\\ University of Manitoba\\
    Winnipeg, Manitoba, Canada}
    \email{charette@cc.umanitoba.ca}
\date{\today}

\begin{abstract}
Given an affine isometry of $\R^3$ with hyperbolic linear part, its
Margulis invariant measures signed Lorentzian displacement along an
invariant spacelike line.  In order for a group generated by hyperbolic
isometries to act properly on $\R^3$, the sign of the Margulis
invariant must be constant over the group.  We show that, in the case
when the linear part is the
fundamental group of a punctured torus, positivity of the
Margulis invariant over any finite generating set does not imply that
the group acts properly.  This contrasts with the case of a pair
of pants, where it suffices to check the sign of the Margulis
invariant for a certain triple of generators.
\end{abstract}

\maketitle 

\section{Introduction}\label{sec:intro}

A flat Lorentz manifold is a quotient $\R^3/\Gamma$, where $\Gamma$ is a group of affine Lorentz transformations acting properly on $\R^3$.   Lorentz transformations preserve a symmetric, non-degenerate bilinear form of signature $(2,1)$, denoted $\ldot{\cdot,\cdot}$.

Margulis~\cite{M1,M2} proved that the fundamental group of a flat Lorentz manifold does not have to be virtually polycyclic, thus answering in the negative a question posed by Milnor.  Milnor had asked whether Auslander's conjecture, that every crystallographic group must be amenable, can be generalized to non-cocompact groups. 

Margulis' examples are affine groups whose linear parts are {\em Schottky groups}: free, discrete groups such that every non-identity element is hyperbolic.  Drumm~\cite{D1} (see also~\cite{CG}) generalized his result by exhibiting fundamental domains for the actions of such groups; he also constructed examples containing parabolics.

While a group of isometries $G$ acts properly on the hyperbolic plane if it is discrete, a group of affine transformations whose linear part is $G$ might not act properly on $\R^3$, even though it acts freely.  We are thus led to study {\em affine deformations} of a discrete group $G$, which are groups of affine transformations whose linear part is $G$.  We seek conditions for an affine deformation of $G$ to act freely and properly on $\R^n$, yielding a flat affine manifold whose fundamental group $\pi$ is isomorphic to $G$.  

We require that the projection of $\pi$ onto its linear part be injective, which amounts to producing an affine deformation of $G$ by assigning translational parts to a set of generators.  Drumm-Goldman~\cite{DG} determined geometric conditions for such an assignment to yield a proper deformation of $G$, when $G$ is a Schottky group.  It leads one to wonder, can properness of an affine deformation of a Schottky group be determined by conditions on some generating set of $G$?

An affine hyperbolic isometry admits a measure of signed Lorentzian displacement, called the {\em Margulis invariant}.  For an affine group of hyperbolic isometries to act properly on $\R^3$, the Margulis invariant must be either always positive in the group, or always negative~\cite{M1,M2}.  

It has been conjectured by Goldman-Margulis that this necessary condition is also sufficient.  An interesting, further question is, if the Margulis invariant is positive for a finite set of elements in the group, can we deduce that the group must act properly on $\R^3$?

Jones~\cite{J} proved that when $G=\langle g,h\rangle$ is the fundamental group of a pair of pants, an affine deformation of the group $\Gamma=\langle \gamma,\eta\rangle$ acts properly and freely if and only if the Margulis invariants of $\gamma$, $\eta$ and $\gamma\eta$ all have the same sign.  This is achieved by showing that the translational parts assigned to $g$ and $h$ satisfy the Drumm-Goldman geometric condition.

In this paper, we show that Jones' result does not extend to the case when $G$ is the fundamental group of a punctured torus.  We elaborate on a construction by Drumm~\cite{D2}, which gives examples of groups $ \langle \gamma,\eta\rangle$, where the Margulis invariant for both $\gamma$ and $\eta$ is positive, but is negative for $\gamma\eta$.  Theorem~\ref{thm:main}, stated below, describes examples of affine groups $\langle \gamma,\eta\rangle$, with Schottky linear part, such that:
\begin{itemize}
\item the Margulis invariant of every word of length up to $n$ is positive;
\item the Margulis invariant of $\gamma\eta^n$ is negative.
\end{itemize}
We can construct such counterexamples for every $n$, with a fixed Schottky linear part.

Recall that every hyperbolic isometry $w$ admits an attracting, future-pointing eigenvector  (of Euclidean length one), which we denote by $\xp{w}$. The fundamental group of a punctured torus has the following interesting property.  If $\langle g,h\rangle$ is the fundamental group of a punctured torus, the set of vectors $\vv$ such that $\ldot{\vv,\xp{gh^i}}>0$, for all $i\leq n$, intersects non-trivially with the set of vectors $\vv$ such that $\ldot{\vv,\xp{gh^{n+1}}}<0$.  This is not true if $\langle g,h\rangle$ is the fundamental group of a pair of pants.

The paper is organized as follows.  Section~\ref{sec:properties}
introduces notation.  We define the Margulis invariant and describe
some of its properties, including a useful addition formula. 

Section~\ref{sec:dynamics} discusses relative positioning of
eigenvectors for words in the Schottky group.  We apply some basic
ideas of symbolic dynamics on the boundary of the hyperbolic plane
(see, for instance~\cite{B}).  We show that there is a non-empty set of vectors $\vv$, such that:
\begin{itemize}
\item $\ldot{\vv,\xo{gh^n}}<0$;
\item $\ldot{\vv,\xo{w}}>0$, for every word $w$ of length up to $n$, except if the left factor of $w$ is of the form $ghg^{-1}$.
\end{itemize}
This ensures that, choosing an appropriate factor of $\vv$ to be the translational part of $\gamma$, the Margulis invariant is always positive, as long as there are no factors of the form $ghg^{-1}$.  

In the second item, the case of words starting with $ghg^{-1}$ poses a
little bit of a difficulty, since the positioning of their
eigenvectors potentially yield negative values of the Margulis
invariant for words of small length.  However, we may impose a certain
technical condition on $\langle g,h\rangle$, called {\em Property C};
roughly speaking, Property C bounds the distance between the
attracting eigenvector of $gh^{-1}$ and the attracting eigenvector of
$ghg^{-1}$.  In the addition formula for the Margulis invariant, negative contributions from words of the form $ghg^{-1}w$ are then cancelled out by contributions from words of the form $gh^{-1}g^{-1}$.  Property C thus ensures that the Margulis invariant is positive for {\em every} word up to length $n$.

Finally, in section~\ref{sec:main}, we prove the main theorem:

\begin{thm}\label{thm:main}
Let $G=\langle g,h\rangle$ be a transversal Schottky group having Property C.  Then for every $n\geq 1$, $G$ admits an affine deformation  $<\gamma,\eta>$ that does not act properly on $\affine$, such that the Margulis invariant of every word of length less than or equal to $n$ is positive.
\end{thm}

\section{Margulis invariant, properties}\label{sec:properties}

Let $\rto$ denote the three-dimensional vector space with the indefinite, symmetric bilinear form defined as follows:
\begin{equation*}
\ldot{\vx,\vy}=x_1y_1+x_2y_2-x_3y_3,
\end{equation*}
where $\vx=(x_1,x_2,x_3)$ and $\vy=(y_1,y_2,y_3)$.  A vector $\vx$ is called {\em spacelike}, {\em timelike} or {\em lightlike}, respectively, if $\ldot{\vx,\vx}$ is positive, negative or null.  

The set of lightlike vectors is called the {\em lightcone}.  The lightcone minus the zero vector admits two connected components, according to the sign of the third coordinate.  A non-zero lightlike vector whose third coordinate is positive (resp. negative) is {\em future-pointing} (resp. {\em past-pointing}).  Denote by $\F$ the set of future-pointing lightlike rays.

The {\em Lorentz-perpendicular plane} of $\vv\in\rto$ is:
\begin{equation*}
\vv^\perp =\{ \vx\in\rto~\mid~\ldot{\vx,\vv}=0\} .
\end{equation*}

A vector $\vv$ is {\em unit-spacelike} if $\ldot{\vv,\vv}=1$.  Denote by $\S$ the set of all unit-spacelike vectors:
\begin{equation*}
\S=\{ \vv\in\rto~\mid~\ldot{\vv,\vv}=1\} .
\end{equation*}

If $\vv$ is spacelike, its Lorentz-perpendicular plane intersects the lightcone in two lightlike lines.

\begin{defn}
Let $\vv\in\S$.  Choose $\xp{\vv},~\xm{\vv}\in\vv^{\perp}$ as follows:
\begin{itemize}
\item $\xp{\vv},~\xm{\vv}\in\F$;
\item $\xp{\vv},~\xm{\vv}$ are of Euclidean length one;
\item $(\vv,\xm{\vv},\xp{\vv})$ is a positively oriented basis.
\end{itemize}
Then $(\vv,\xm{\vv},\xp{\vv})$ is called the {\em null frame} relative to $\vv$.
\end{defn}
(Requiring that $\xpm{\vv}$ are unit length simply ensures uniqueness of choice.)

The {\em Lorentz cross-product} is the unique bilinear transformation $\boxtimes:\rto\times\rto\longrightarrow\rto$ such that:
\begin{equation*}
\ldot{\vu,\vv\boxtimes\vw}=\det[\vu~\vv~\vw] .
\end{equation*}
In particular, $\ldot{\vv,\vv\boxtimes\vw}=0$.  It is easily seen that $\vv\in\S$ is a positive scalar multiple of $\xm{\vv}\boxtimes\xp{\vv}$.

We may identify $\rto$ with $\affine$, the three-dimensional affine space whose space of translations is $\rto$.

\subsection{Isometries}

Let $\isoaff$ be the group of affine isometries of $\affine$.  Via the identification between $\rto$ and $\affine$, any $\gamma\in\isoaff$ may be written as follows:
\begin{equation*}
\gamma (x)=g(x)+\vg,
\end{equation*}
where $g\in\oto$ and $\vg\in\rto$; $g$ is called the {\em linear part} of $\gamma$ and $\vg$, its {\em translational part}.  

We shall call $\gamma$ an {\em affine deformation of $g$}.  More generally, $\Gamma\subset\isoaff$ is called an {\em affine deformation} of $G$ if its linear part is $G$.  An affine deformation is obtained from a linear group by assigning translational parts to its generators.

Isometries in $\oto$ either preserve each connected component of the lightcone or interchange them.  Let $\soto$ denote the connected component of the identity in $\oto$; $\soto$ consists of orientation- and time orientation-preserving isometries.  

Let $\gamma\in\isoaff$; suppose its linear part $g$ lies in $\soto$.  Both $\gamma$ and $g$ are called {\em hyperbolic} if $g$ admits three distinct eigenvalues.  Denote its smallest eigenvalue by $\lag$; it must be positive and less than 1 and the eigenvalues of $g$ are thus $\lag,1,\lag^{-1}$.  The $\lag^{\pm 1}$ eigenvectors must be lightlike and the 1-eigenvectors, spacelike.  Choose $\xp{g}$, $\xm{g}$ and $\xo{g}$ as follows:
\begin{itemize}
\item $\xm{g}$ is a $\lag$-eigenvector and $\xp{g}$, a $\lag^{-1}$-eigenvector;
\item $\xo{g}\in\S$ is an eigenvector such that $\xpm{\xo{g}}=\xpm{g}$.
\end{itemize}
This is the null frame of $\xo{g}$ and is uniquely determined by $g$.  Note that $\xo{g^{-1}}=-\xo{g}$ and that for any isometry $h\in\soto$:
\begin{equation*}
\xo{hgh^{-1}}=h\xo{g}.
\end{equation*}

\subsection{The Margulis invariant}

Suppose $\gamma\in\isoaff$ is hyperbolic, with linear part $g$.  Then $\gamma$ preserves a unique line in $\affine$ that is parallel to $\xo(g)$.  (It is the only $\gamma$-invariant line if $\gamma$ acts freely.)  Since $(\xm{g},\xp{g},\xo{g})$ forms a basis for $\rto$, the following definition is independent of the choice of $p$.
\begin{defn}
Let $\gamma\in\isoaff$ be hyperbolic with linear part $g$.  The {\em Margulis invariant} of $\gamma$ is:
\begin{equation}
\alpha(\gamma)=\ldot{\gamma(p)-p,\xo{g}},
\end{equation}
where $p$ is an arbitrary point in $\affine$.
\end{defn}
In particular, if $\vg$ is the translational part of $\gamma$, then $\alpha(\gamma)=\ldot{\vg,\xo{g}}$.
This invariant has the following properties: 

\begin{itemize}
\item $\alpha(\gamma^{-1})=\alpha(\gamma)$;
\item For any $\eta\in\isoaff$, $\alpha(\eta\gamma\eta^{-1})=\alpha(\gamma)$.
\end{itemize}

The Margulis invariant yields a useful criterion for determining that an action is non-proper.

\begin{lemma} (Margulis~\cite{M1,M2})\label{alphasign}
Let $\gamma$, $\eta\in\isoaff$ be hyperbolic transformations such that $\alpha(\gamma)\alpha(\eta)<0$.  Then $\langle\gamma,\eta\rangle$ does not act properly on $\affine$.
\end{lemma}

\subsection{The Margulis invariant of a cyclically reduced word}
Suppose $g=\langle g_1,\ldots,g_n\rangle$ is freely generated by the $g_i$'s.  Then every $g\in G$ can be uniquely written as a {\em word} in the generators:
\begin{equation*}
w=g_{i_k}^{j_k}g_{i_{k-1}}^{j_{k-1}}\ldots g_{i_1}^{j_1},
\end{equation*}
where $1\leq i_l\leq n$, $j_l=\pm 1$ and $g_{i_l}\neq g_{i_{l+1}}^{-1}$.  This last condition states that $w$ is a {\em reduced word}.  The length of a reduced word is well-defined.  We call $g_{i_k}$ the {\em terminal letter} of $w$ and $g_{i_1}$, its {\em initial letter}.  (This is consistent with a left action of $\isoaff$ on $\E$.)

We will say that a reduced word $g$ is {\em cyclically reduced} if $g^2$ is reduced.  In other words, $g$ is not the conjugate of a reduced word of shorter length. 

\begin{lemma} (Drumm-Goldman~\cite{DG})\label{lem:alphasum}
Let $\Gamma=\langle\gamma_1,\ldots,\gamma_n\rangle$ be an affine
deformation of a group that is freely generated by hyperbolic elements.  Denote the linear part of $\gamma_i$ by $g_i$ and its translational part by $\vv_{g_i}$.  Suppose $\gamma\in\Gamma$ is a cyclically reduced hyperbolic isometry and write its linear part as the reduced word:
\begin{equation*}
w=g_{i_1}^{j_1}g_{i_2}^{j_2}\ldots g_{i_m}^{j_m},
\end{equation*}
where $j_i=\pm 1$. Then:
\begin{equation}
\alpha(\gamma)=\sum_{k=1}^n\ldot{\vv_{g_{i_k}},\xo{w_k}},
\end{equation}
where $w_k$ is a cyclically reduced word whose terminal letter is $g_{i_k}$.
\end{lemma}
Note that we have written the terminal letter as $g_{i_1}^{j_1}$, in order to simplify notation in the proof.
\begin{proof}
To simplify notation, set:
\begin{align*}
h_k &=g_{i_{k}}^{j_{k}}\\
\vv_k&=\vv_{h_k},
\end{align*}
for $1\leq k\leq m$.  Thus:
\begin{align*}
\alpha(\gamma) &= \ldot{\vv_1,\xo{w}}+\sum_{k=2}^n\ldot{h_1\ldots h_{k-1}(\vv_{k}),\xo{w}} \\
&=   \ldot{\vv_1,\xo{w}}+\sum_{k=2}^n\ldot{\vv_{k},(h_1\ldots h_{k-1})^{-1}\xo{w}} \\
&= \ldot{\vv_1,\xo{w}}+\sum_{k=2}^n\ldot{\vv_k,\xo{h_k\ldots h_mh_1\ldots h_{k-1}}},
\end{align*}
since $w_2^{-1}\xo{w_1}=\xo{w_2^{-1}w_1w_2}$.

Now consider each summand $\ldot{\vv_k,\xo{h_k\ldots h_mh_1\ldots h_{k-1}}}$.  If $j_{k}=1$, then $h_k=g_{i_k}$ and $\vv_k=\vv_{g_{i_k}}$ so that the summand is of the desired form.  Moreover, $w_k$ is reduced because $\gamma$ is cyclically reduced, and cyclically reduced because $\gamma$ is reduced.

If $j_{k}=-1$, then the corresponding summand is:
\begin{align*}
\ldot{-g_{i_k}^{-1}(\vv_{g_{i_k}}),\xo{g_{i_k}^{-1}\ldots h_mh_1\ldots h_{k-1}}} 
&= -\ldot{\vv_{g_{i_k}},\xo{w_k'g_{i_k}^{-1}}} \\
&= \ldot{\vv_{g_k},\xo{g_{i_k}(w_k')^{-1}}} ,
\end{align*}
where $w_k'=h_{k+1}\ldots h_mh_1\ldots h_{k-1}$.  As above, $w_k=g_{i_k}(w_k')^{-1}$ is both reduced and cyclically reduced.
\end{proof}
\begin{ex}
Let $\gamma$, $\eta\in\isoaff$ be two hyperbolic elements with respective linear parts $g$, $h$ and respective translational parts $\vg$, $\vh$.  Then:
\begin{align*}
\alpha(\gamma\eta) &=\ldot{\vg,\xo{gh}}+\ldot{\vh,\xo{hg}} \\
\alpha(\gamma\eta^{-1}) &=\ldot{\vg,\xo{gh^{-1}}}+\ldot{\vh,\xo{hg^{-1}}} .
\end{align*}
More generally, for $n\geq 1$: 
\begin{equation}
\alpha(\gamma\eta^n)=\ldot{\vg,\xo{gh^n}}+\sum_{k=0}^{n-1}\ldot{\vh,\xo{h^{n-k}gh^k}}.
\end{equation}
\end{ex}
Following Drumm~\cite{D2}, we can use the formula for
$\alpha(\gamma\eta)$ to
construct non-proper actions.  Indeed, the intersection of the half-spaces: 
\begin{equation*}
\{ \vv\in\rto\mid\ldot{\vv,\xo{g}}>0\}\bigcap\{ \vv\in\rto\mid
\ldot{\vv,\xo{gh}}<0\}
\end{equation*}
is convex and non-empty.  Choose any $\vv$ in this
intersection and choose the translational part of $\eta$ to be any
$\vh$ such that $\alpha(\eta)$ is positive.  Then set $\vg =k\vv$,
where $k$ is large enough so that $\alpha(\gamma\eta)<0$.  Observe
that $\alpha(\gamma)$ and $\alpha(\eta)$ are both positive, but by
Lemma~\ref{alphasign}, the group generated by $\gamma$ and $\eta$ does
not act properly, since the subgroup $\langle\gamma ,\gamma\eta\rangle$ does not act properly.

\section{Symbolic dynamics on the boundary at infinity}\label{sec:dynamics}

We now discuss the relative positioning of attracting null eigenvectors of words, based on the positioning of the eigenvectors of a generating set for the group.  In the case of a {\em Schottky group}, defined below, we can use Brouwer's fixed point theorem.

Denote the {\em sign} of an integer $i$ by $\sigma(i)$ and the closure of a set $A$ by $cl(A)$. In what follows, the closure of subsets of $\F$, the set of future-pointing lightlike rays, will always be taken relative to $\F$. 

\begin{defn}
A connected set of future-pointing lightlike rays is called a {\em conical interval}.  A {\em conical neighborhood} of $\vx\in\F$ is a conical interval containing $\vx$.
\end{defn}

\begin{figure}
\input{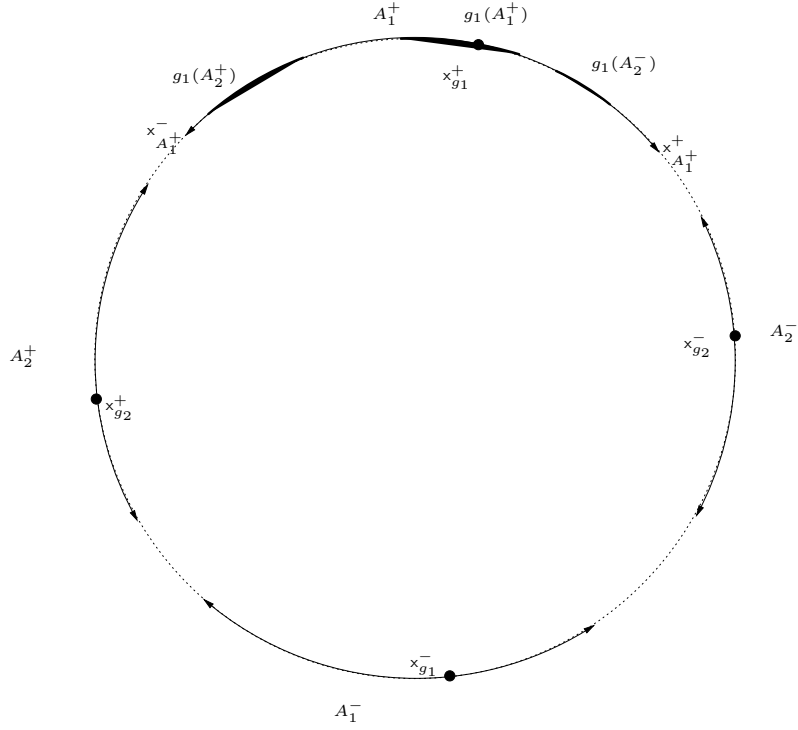}
\caption{A Schottky system for $\langle g_1,g_2\rangle$.  Since $g_1(B)\subset A^+_1$ for $B=A^+_1,~A^\pm_2$, $\xp{h}\in A^+_1$ for $h=g_1^2,~g_1g_2,~g_1g_2^{-1}$.}
\label{fig:Schottky}
\end{figure}

\begin{defn}
Let $g_1,\ldots,g_n\in\soto$.  Then $\langle g_1,\ldots,g_n\rangle$ is a {\em Schottky group} if there exist $2n$ disjoint closed conical intervals $A^{\pm}_i$ such that, for $1\leq i\leq n$:
\begin{equation*}
g_i(A_i^-)=cl(\F - A_i^+).
\end{equation*}
The generators $g_1,\ldots,g_n$ are called {\em Schottky generators} of $G$ and the conical intervals $A^\pm_i$ are called a {\em Schottky system}.
\end{defn}
See Figure~\ref{fig:Schottky}.

Note that $\xpm{g_i}\in A_i^\pm$.  In fact, as is illustrated in Figure~\ref{fig:Schottky}, Brouwer's fixed point theorem implies that when $w=g_{i_n}^{j_n}\ldots g_{i_1}^{j_1}$ is a cyclically reduced word, then $\xp{w}\in A_{i_n}^{\sigma(j_n)}$ and $\xm{w}\in A_{i_1}^{\sigma(-j_1)}$.

Equivalently, a finitely generated subgroup of $\oto$ is Schottky if and only if it is free, discrete and purely hyperbolic.

Let $U\neq\emptyset,~\F$ be a closed conical interval.  Denote by $\xm{U}$, $\xp{U}$, respectively, the (Euclidean) unit-length lightlike vectors spanning the boundary of $U$, such that $\ldot{\vx,\xm{U}\boxtimes\xp{U}}>0$ for $\vx\in U$.

If $\UM$, $\UP$ is a disjoint pair of closed and non-trivial conical intervals, define:
\begin{align*}
\C{\UM,\UP}=\{\vv\in\rto~\mid~ & \ldot{\vv,\vu}>0\mbox{~for all~}\vu\in\S~\mbox{such that~}\\ 
                               & \xm{\vu}\in\UM,~\xp{\vu}\in\UP\} .
\end{align*}

This is a non-empty intersection of half-spaces.  In fact, it is simply the intersection of four half-spaces, namely, those determined by $\xpm{\UM}\boxtimes\xpm{\UP}$.

\begin{lemma}
Let $\UM$, $\UP$ be closed, disjoint, non-trivial conical intervals.  Then $\C{\UM,\UP}$ is the cone of positive linear combinations of $\xp{\UM}$, $\xm{\UP}$, $-\xp{\UP}$ and $-\xm{\UM}$.
\end{lemma}
\cqfd 

See Figure~\ref{fig:cone}.  In particular, if $W$ is a set of words $g$ such that $\xm{g}\in\UM$ and $\xp{g}\in\UP$, then every $\vv\in\C{\UM,\UP}$ satisfies $\ldot{\vv,\xo{g}}>0$, for every $g\in W$.

\begin{figure}
\input{cone.pstex_t}
\caption{The cone $\C{\UM,\UP}$, which is the set of positive linear combinations of $\xm{\UP}$, $\xp{\UM}$, $-\xm{\UM}$ and $-\xp{\UP}$.  For any $\vu$ such that $\xm{\vu}\in\UM$ and $\xp{\vu}\in\UP$, and any $\vv\in\C{\UM,\UP}$, $\ldot{\vu,\vv}>0$.}
\label{fig:cone}
\end{figure}

\begin{ex}
Let $G=\langle g,h\rangle$ be a Schottky group with Schottky system $A^\pm_g$, $A^\pm_h$.  Set $\UP=A^+_h$ and let $\UM$ be the smallest (closed) conical interval containing $A^-_h\bigcup A^\pm_g$.  If $\vv\in\C{\UM,\UP}$:
\begin{equation*}
\ldot{\vv,\xo{w}}>0,
\end{equation*}
for every word $w$ with terminal letter $h$.
Drumm-Goldman~\cite{DG} construct affine deformations of Schottky groups that act properly on $\affine$, using this idea.
\end{ex}

\subsection{Ordering on the future lightcone}

Given a Schottky subgroup of $\soto$, we seek an ordering of the
attracting eigenvectors in the Schottky intervals.  Let 
$\vx,\vy,\vz\in\F$; we write:
\begin{equation*}
\PO{\vx,\vy,\vz }
\end{equation*}
if and only if $(\vx,\vy,\vz)$ is a right-handed positively-oriented
basis, i.e. $\det([\vx~\vy~\vz])>0$.  More generally, $\PO{\vv_1,\ldots,\vv_k}$ if and only if $\PO{\vv_{i_1},\vv_{i_2},\vv_{i_3}}$ for every triple $1\leq i_1<i_2<i_3\leq k$.  

We extend the notation to conical intervals in the obvious manner.  Thus, if $U$ is a conical interval, we will write $\PO{\vx,U,\vz}$ to mean that $\PO{\vx,\vy,\vz}$ for every $\vy\in U$.

Obviously, $\PO{\vx,\vy,\vz }$ if and only if $\PO{\vy,\vz,\vx }$.

\begin{lemma}
Let $\vx,\vy,\vz\in\F$.  Then $\PO{\vx,\vy,\vz }$ if and only if $\ldot{\vx,\vy\boxtimes\vz}>0$.
\end{lemma}
\cqfd 

If $g\in\soto$, then $g$ is orientation-preserving; thus $\PO{\vx,\vy,\vz }$ if and only if $\PO{g(\vx),g(\vy),g(\vz) }$.

We will order vectors in a conical interval according to one of its boundary components.  Note that if $A$ is a conical interval and $\vx$, $\vy\in A$,

\begin{equation*}
\PO{\xp{A},\vx,\vy}\Longleftrightarrow\PO{\xm{A},\vx,\vy},
\end{equation*}
which both mean that, starting at $\xp{A}$ and moving
counterclockwise, we meet $\vx$, $\vy$ and then $\xm{A}$.

Recall that if $G=\langle g_1,\ldots,g_n\rangle$ is a Schottky group,
$\xp{w}\in A^\pm _i$ for every reduced word $w$ with terminal letter $g_i^{\pm
  1}$.  The following lemma yields an ordering of the attracting
eigenvectors of cyclically reduced words of the form $h^{\pm
  1}g_{i_1}^{\pm 1}w_1$ and $hg_{i_2}^{\pm 1}w_2$, where $1\leq
i_1,i_2\leq n$ and $j_1,j_2$ are non-zero integers.  For simplicity,
we write $g_1=g_{i_1}^{\pm 1}$ and  $g_2=g_{i_2}^{\pm 1}$.  This is
allowed, since any Schottky generator can be replaced by its inverse
and the generators can be reindexed.

\begin{lemma}\label{lem:order}
Let $G=\langle g_1,\ldots,g_n\rangle$ be a Schottky group and let
$h=g_i^k$, for some $1<i<n$ and $k=\pm 1$.  Set $A=A^{\sigma(k)}_i$.
Let $w_1$, $w_2$ be non-trivial words in $G$ such that $hg_1w_1$ and $hg_2w_2$ are cyclically reduced.    Then:
\begin{equation*}
\PO{\xp{A},\xp{hg_1w_1},\xp{hg_2w_2}}\Longleftrightarrow\PO{\xm{h},\xp{g_1},\xp{g_2}}.
\end{equation*}
\end{lemma}

\begin{proof}
Set $B_j=g_j(w_j(A))$, $j=1,2$ and $A^-= A_i^{\sigma(-k)}$.  Then $B_j\subset A^+_j$, which is disjoint from $A^-$.  Thus $\PO{\xm{h},B_1,B_2}$ makes sense and:
\begin{align*}
\PO{\xm{h},\xp{g_1},\xp{g_2}}&\Longleftrightarrow\PO{\xm{h},B_1,B_2} \\
& \Longleftrightarrow\PO{\xm{A^-},B_1,B_2} \\
& \Longleftrightarrow\PO{\xp{A},h(B_1),h(B_2)},
\end{align*}
since $h(\xm{A^-})$ is parallel to $\xp{A}$.  By Brouwer's fixed point theorem, $h(B_i)$ contains the attracting eigenvector of $hg_iw_i$ and the result follows.
\end{proof}

\begin{cor}\label{cor:order}
Using the hypotheses of Lemma~\ref{lem:order}, let $w\in G$ be a word with terminal letter $g=g_i^k$, such that $w_1'=whg_1w_1$ and $w_2'=whg_2w_2$ are cyclically reduced.  Then: 
\begin{equation*}
\PO{\xp{A'},\xp{w_1'},\xp{w_2'}}\Longleftrightarrow \PO{\xm{h},\xp{g_1},\xp{g_2}},
\end{equation*}
where $A'=A^{\sigma(k)}_i$.
\end{cor}
\begin{proof}
This follows from Lemma~\ref{lem:order} and the observation that $\PO{\xm{A''},\xp{hg_1w_1},\xp{hg_2w_2}}$ if and only if $\PO{\xp{A},\xp{hg_1w_1},\xp{hg_2w_2}}$, where $A''=A^{\sigma(-k)}_i$.
\end{proof}

\subsection{Transversal isometries}

\begin{defn}
Suppose $<g,h>\subset\soto$ is a Schottky group.  The group is said to be {\em transversal} if the planes ${\xo{g}}^\perp$ and ${\xo{h}}^\perp$ intersect in a timelike line.
\end{defn}
Thinking of the interior of the projectivized lightcone as the hyperbolic plane, this means that $<g,h>$ is the fundamental group of a punctured torus.

For the remainder of this section, we will assume $\langle g,h\rangle$ to be transversal.  Taking $h^{-1}$ if necessary, we may further assume that:
\begin{equation*}
\PO{\xm{g},\xm{h},\xp{g}}.
\end{equation*}
Moving counterclockwise on the future lightcone and starting at $\xm{g}$, we meet, in order: $\xm{h}$, $\xp{g}$ and finally, $\xp{h}$.

\begin{lemma}\label{lem:transv order}
Let $\langle g,h\rangle$ be a transversal Schottky group with Schottky system 
$\{ A^\pm_g,A^\pm_h\}$, with $\PO{\xm{g},\xm{h},\xp{g}}$.  Set $\vx^+=\xp{A^+_g}$.  Then the following hold:
\begin{enumerate}
\item $\PO{\vx^+,\xp{gh^{-1}w_1},\xp{g^2w_2},\xp{ghw_3}}$ for every $w_1$, $w_2$, $w_3$ yielding cyclically reduced words;
\item $\PO{\vx^+,\xp{gh^igw},\xp{gh^{n+1}}}$, $i\leq n$, for any word $w$ such that $gh^igw$ is cyclically reduced; 
\item $\PO{\vx^+,\xp{gh^i},\xp{gh^{i+j}}}$, for all $i\geq 0$ and $j>0$.
\end{enumerate}
\end{lemma}

\begin{figure}
\input{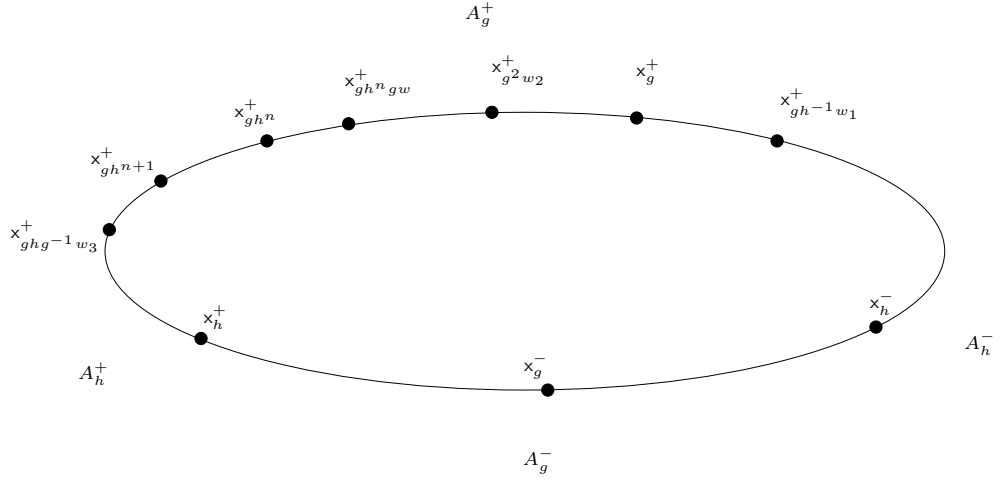}
\caption{The ordering of some attracting eigenvectors in the case of
  transversal Schottky generators $g$ and $h$.}
\label{fig:order}
\end{figure}

See Figure~\ref{fig:order}.

\begin{proof}
 Items 1 and 2 follow from Lemma~\ref{lem:order}.

To prove Item 3, consider the conical interval $U$ contained in $A^+_g$ bounded by $\xp{gh^i}$ and $\xm{A^+_g}$.  Then $h^j(U)\subset A^+_h$ and thus, $\PO{\xp{A^-_g},h^j(U),\xp{gh^i}}$.  It follows that $\PO{\vx^+,gh^{i+j}(U),\xp{gh^i}}$.  By Brouwer's fixed point theorem, $\xp{gh^{i+j}}\in gh^{i+j}(U)$.
\end{proof} 

Let $\UM$ be the smallest conical interval containing $A^\pm_h$ and $A^-_g$.  Next, for each $n\geq 1$, set $\UPN$ to be the conical interval contained in $A^+_g$ bounded by $\xp{A^+_g}$ and $\xp{gh^{n-1}}$.  Thus:
\begin{align*}
\xm{\UM}  & =\xm{A^-_h}       & \xp{\UM}      & =\xp{A^+_h}\\
\xm{\UPN} & =\xp{gh^{n-1}}    & \xp{\UPN} & =\xp{A^+_g}.
\end{align*}

By Lemma~\ref{lem:transv order}, if $\vv\in\C{\UM,\UPN}$, then $\ldot{\vv,\xo{g}}>0$ and, furthermore, $\ldot{\vv,\xo{w}}>0$ for every cyclically reduced word $w$ of length less than or equal to $n+1$, with terminal factor $g^k$, $k\geq 2$, $gh^ig$ or $gh^{-1}$.

Now, for each $n\geq 1$, set:
\begin{equation*}
\HN =\{ \vv\in\rto~\mid~\ldot{\vv,\xo{gh^n}}<0\}.
\end{equation*}

\begin{figure}
\input{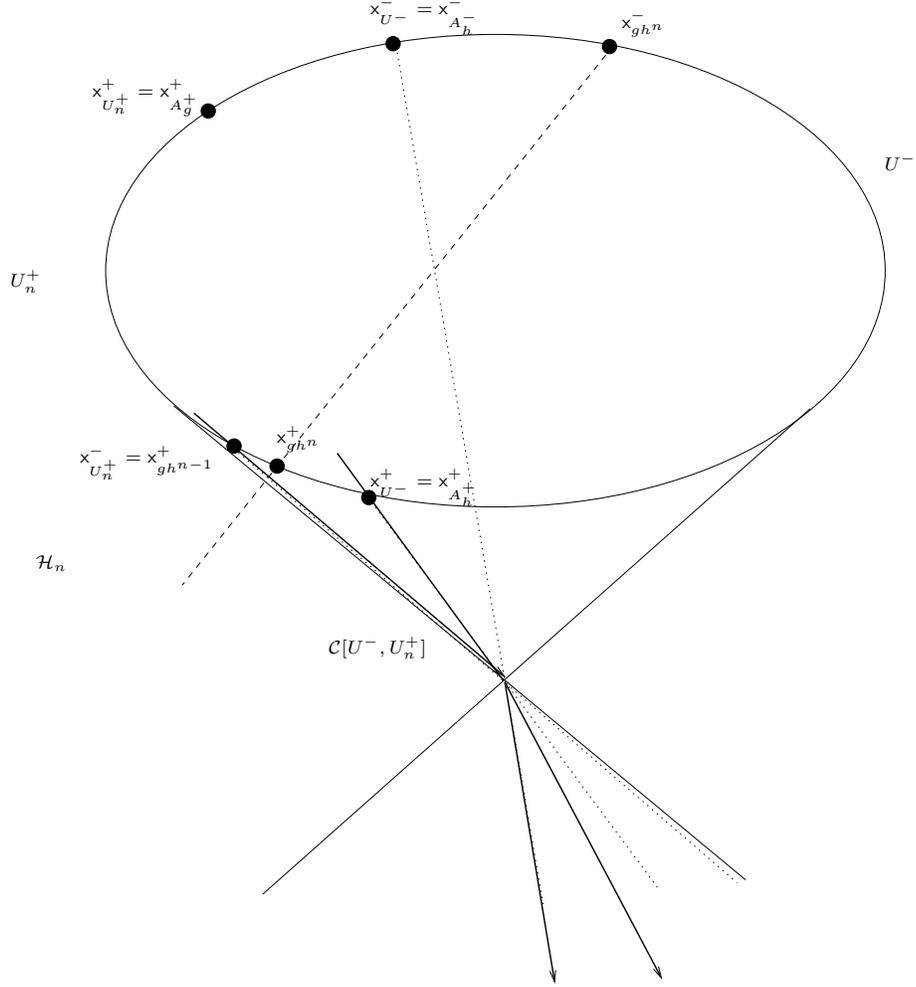}
\caption{The cone $\C{\UM,\UPN}$, bounded by the darker lines, intersects non-trivially with the half-space $\HN$, bounded by the dashed line.}
\label{fig:cone2}
\end{figure}

Observe that $\HN\bigcap\C{\UM,\UPN}$ is non-empty, since both sets contain the conical interval in $A^+_g$ bounded by $\xp{gh^{n-1}}$ and $\xp{gh^{n}}$.  (See Figure~\ref{fig:cone2}.)

\subsection{Words of the form $gh^ig^{-1}w$}
Our goal is to assign translational parts $\vg$, $\vh$ to $g,~h$
respectively, to obtain an affine deformation of a Schottky group
$\Gamma=\langle\gamma,\eta\rangle$, such that 
\begin{itemize}
\item $\alpha(\omega)>0$ for every $\omega\in\Gamma$ of length less
  than or equal to $n$, and 
\item $\alpha(\gamma\eta^n)<0$.
\end{itemize}
 By Lemma~\ref{lem:alphasum}:
\begin{equation}\label{eq:alphasum}
\alpha(\omega)=\sum_{k=1}^n\ldot{\vv_k,\xo{w_k}},
\end{equation}
where $\vv_k=\vg$ (resp. $\vv_k=\vh$) and $w_k$ is a cyclically reduced word whose terminal letter is $g$ (resp. $h$).

We can start by choosing $\vh$ such that $\ldot{\vh,\xo{w}}>0$ for every cyclically reduced word $w$ with terminal letter $h$.  Then, we can choose any $\vv\in\HN\bigcap\C{\UM,\UPN}$, then set $\vg$ to be a scalar multiple of $\vv$ satisfying $\alpha(\gamma\eta^n)<0$.

If every summand in Equation~\eqref{eq:alphasum} were positive when the
length of $\gamma$ is at most $n$, then we would be done.  And
$\ldot{\vg,\xo{w}}$ is indeed positive for many choices of $w$.  But
our choice of $\vg$ yields a negative value for $\ldot{\vg,\xo{gh^ig^{-1}}}$.  As a matter of fact,
for $n$ large enough, it is impossible to force all the
summands in Equation~\eqref{eq:alphasum} to be positive.

However, it is possible that every negative summand is cancelled out by a positive summand, in light of the following observation. 

\begin{lemma}\label{lem:conjugate}
Write $\alpha(\omega)$ as in Equation~\eqref{eq:alphasum}.  Then summands of the form $\ldot{\vg,\xo{gh^ig^{-1}w'}}$ are in one-to-one correspondence with summands of the form $\ldot{\vg,\xo{gh^{-i}g^{-1}w''}}$.
\end{lemma}

\begin{proof}
Let $w$ be the word in $g,~h,~g^{-1},~h^{-1}$ corresponding to the linear part of $\gamma$.  As is demonstrated in the proof of Lemma~\ref{lem:alphasum}, a summand of the form $\ldot{\vg,\xo{gh^ig^{-1}w'}}$ may arise from a $g$ factor or a $g^{-1}$ factor.

{\bf Case 1:~}$w=w_1gw_2$, where $g$ induces a summand of the form $\ldot{\vg,\xo{gh^ig^{-1}w'}}$.  Then $w=w_1gh^ig^{-1}w_2'$.  Thus there is also a summand $\ldot{w_1gh^i(-g^{-1}(\vg),\xo{w}}=\ldot{\vg,\xo{gh^{-1}g^{-1}w''}}$.

{\bf Case 2:~}$w=w_1g^{-1}w_2$, where $g^{-1}$ induces a summand of the form $\ldot{\vg,\xo{gh^ig^{-1}w'}}$.  This happens when $w_1=w_1'gh^{-i}$.  Thus there is also a summand $\ldot{w_1'(\vg),\xo{w}}=\ldot{\vg,\xo{gh^{-i}g^{-1}w_2(w_1')^{-1}}}$.
\end{proof}

Thus every negative summand is potentially cancelled out by a positive summand.

{\bf Property C.}  Let $W^+$ (resp. $W^-$) denote the set of all cyclically reduced words of the form $gh^ig^{-1}w$ (resp. $gh^{-i}g^{-1}w$).  Let $A$ be the conical interval in $A^+_g$ bounded by $\xp{gh}$ and $g(\xp{h})$.  
We will say that the transversal Schottky group $\langle g,h\rangle$ satisfies Property C if for every $\vx\in A$: 
\begin{equation*}
-\ldot{\vx,\xo{w_1}}<\ldot{\vx,\xo{w_2}},~\mbox{for all~}w_1\in W^+, w_2\in W^- .
\end{equation*}

Note that $A$ contains $\xp{gh^i}$ for all $i>0$.   Furthermore, we have chosen $A$ so that both terms in the previous inequality are positive.  

Property C is more restrictive than is really necessary, but it will make arguments easier.  Groups that satisfy Property C are easy find: for instance, setting $\xo{g}=(1,0,0)$ and $\lag=e^{-1}$, and $\xo{h}=(0,1,0)$ and $\lah=e^{-2}$, the group $\langle g,h\rangle$ has Property C.

\section{Main theorem}\label{sec:main}

Recall the statement of Theorem~\ref{thm:main}.

\begin{thm}
Let $G=\langle g,h\rangle$ be a transversal Schottky group having Property C.  Then for every $n\geq 1$, $G$ admits an affine deformation  $<\gamma,\eta>$ that does not act properly on $\affine$, such that the Margulis invariant of every word of length less than or equal to $n$ is positive.
\end{thm}

\begin{proof}
Let $n\geq 1$.  By choosing appropriate translational parts $\vg$ and $\vh$, we will exhibit a group $<\gamma,\eta>$, as in the statement of the theorem, such that $\alpha(\gamma\eta^{n})<0$.  

First, we choose a translational part $\vh$ for $h$, such that $\ldot{\vh,\xo{w}}>0$, for every cyclically reduced word $w$ with terminal letter $h$.

The set $A\bigcap\HN\bigcap\C{\UM,\UPN}$ is non-empty, since it contains the conical interval in $A^+_g$ bounded by $\xp{gh^{n-1}}$ and $\xp{gh^n}$.

Choose any $\vv\in A\bigcap\HN\bigcap\C{\UM,\UPN}$.  Lemma~\ref{lem:order} and Property C ensure that,  for any $k>0$, choosing $\vg=k\vv$ will yield $\alpha(\gamma)>0$ for any word of length up to $n$.  But taking $k$ large enough so that $\alpha(\gamma\eta^n)<0$, it follows that $\langle\gamma,\eta\rangle$ cannot act properly on $\affine$.
\end{proof}

\makeatletter \renewcommand{\@biblabel}[1]{\hfill#1.}\makeatother


\begin{thebibliography}{9}

\bibitem{B} ``Ergodic Theory, Symbolic Dynamics and Hyperbolic Spaces,\ '' 
            Bedford, T., Keane, M. and Series, C. eds, Oxford (1991),
            369 pages.
\bibitem{CG}
    Charette, V.\ and Goldman, W.,
    {\em Affine Schottky groups and crooked tilings,\/}
    Proceedings of the Workshop on ``Crystallographic Groups and
        their Generalizations II,\ ''Contemp.\ Math.\ (to appear).

\bibitem{D1}
    Drumm, T.,
    {\em Fundamental polyhedra for Margulis space-times,\/}
    Topology {\bf 31} (4) (1992), 677-683.

\bibitem{D2}
    \bysame,
    {\em Examples of nonproper affine actions,\/}
    Mich.\ Math.\ J.\ {\bf 39} (1992), 435--442.

\bibitem{DG}
    Drumm, T.\  and Goldman, W., 

\bibitem{J}
    Jones, C.,
    doctoral dissertation, University of Maryland, 2003.
    
\bibitem{M1}
    Margulis, G.,
    {\em Free properly discontinuous groups of affine transformations,\/}
        Dokl.\ Akad.\ Nauk SSSR {\bf 272} (1983), 937--940.

\bibitem{M2}
    \bysame,
    {\em Complete affine locally flat manifolds with a free
    fundamental group,\/}
    J.\ Soviet Math. {\bf 134} (1987), 129--134.


\end{thebibliography}
\end{document}